\documentclass[12pt]{amsart}
\usepackage{tikz}
\usepackage{amsmath}
\usepackage{amsthm}
\usepackage{amssymb}
\usepackage{amscd}
\usepackage{bbold}
\usepackage{caption}

\usepackage{enumerate}
\usepackage{graphicx}
\usepackage{fancyhdr}
\usepackage{multicol}
\usepackage[margin={1cm,1cm}]{geometry}

\usepackage{subcaption}
\usepackage{nameref}
\usepackage{algpseudocode,algorithm,algorithmicx}

\usepackage [english]{babel}
\usepackage [autostyle, english = american]{csquotes}
\MakeOuterQuote{"}

\title[Wind tree model]
{Wind-tree model for billiard motion from a signal processing viewpoint}

\date{}

\subjclass[]{}

\author{Enrico Au-Yeung}
\address{Department of Mathematical Sciences\\
          DePaul University\\
         Chicago, IL 60614\\
         USA}
 \email{eauyeun1@depaul.edu}

\author{Nick Kreissler}
\email{nick@hmcproducts.com}

\begin{document}

\begin{abstract}

In the Ehrenfest wind tree model, a point  particle moves on the plane and collides  with randomly placed fixed square obstacles under the usual law of geometric optics.  The particle represents the wind and the squares are the trees. We examine the periodic version of the model. Previous authors analyze the dynamical properties of the model using techniques from algebraic topology or ergodic theory.  In contrast to these works, we adopt a signal processing viewpoint.   We describe the phenomenon of the long-term trajectories by using a 3-state hidden Markov model.

\end{abstract}

\maketitle

		
		\section{Introduction and Motivation}
			
			In the  wind tree model introduced by P. Erhenfest and T. Ehrenfest in 1912, a point particle moves on the plane and collides  with randomly placed fixed square obstacles under the usual law of geometric optics. The particle represents the wind and the squares are the trees.  We investigate the periodic version of the model, where  identical square obstacles are placed periodically in the plane.  Our aim is to understand some of
it dynamical  properties (see Figure 1 for the position of the billiard after 15 collisions with the obstacles.)

\begin{figure}[h!]
   \centering
   \begin{subfigure}[b]{0.4\linewidth}
       \includegraphics[width=\linewidth]{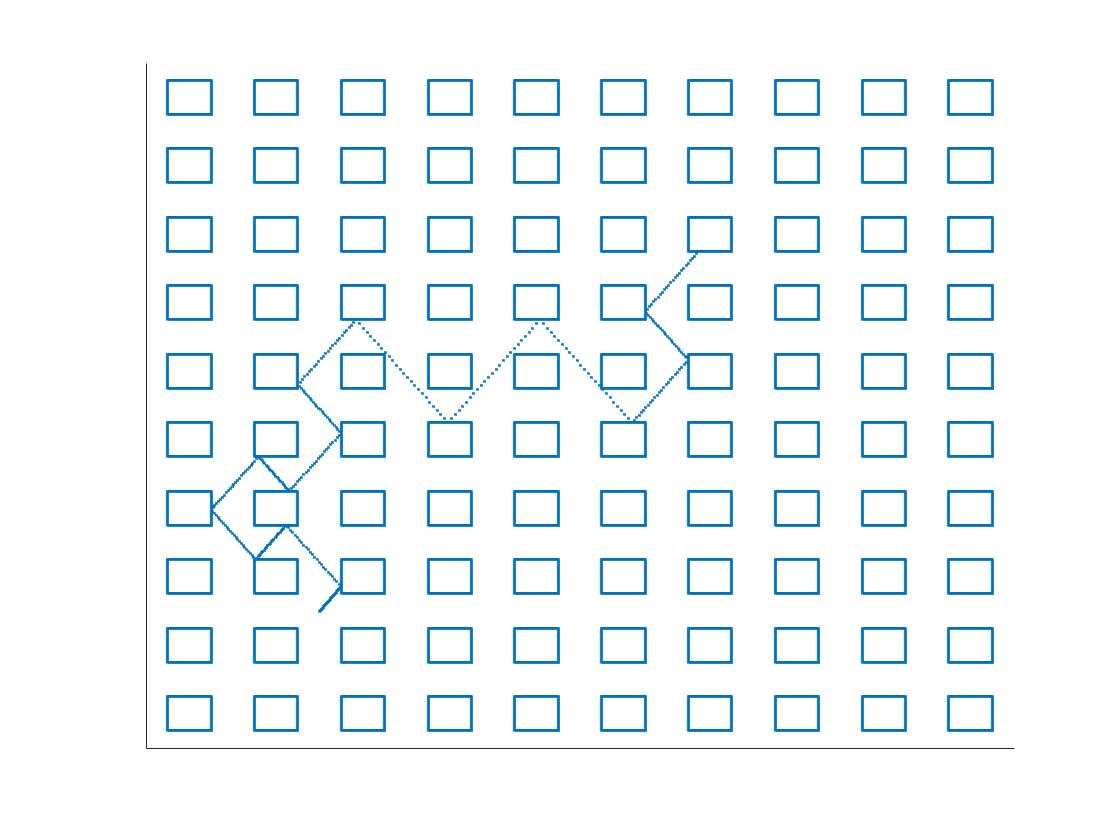}
       \caption{initial slope = 1.414}
   \end{subfigure}
   \begin{subfigure}[b]{0.4\linewidth}
       \includegraphics[width=\linewidth]{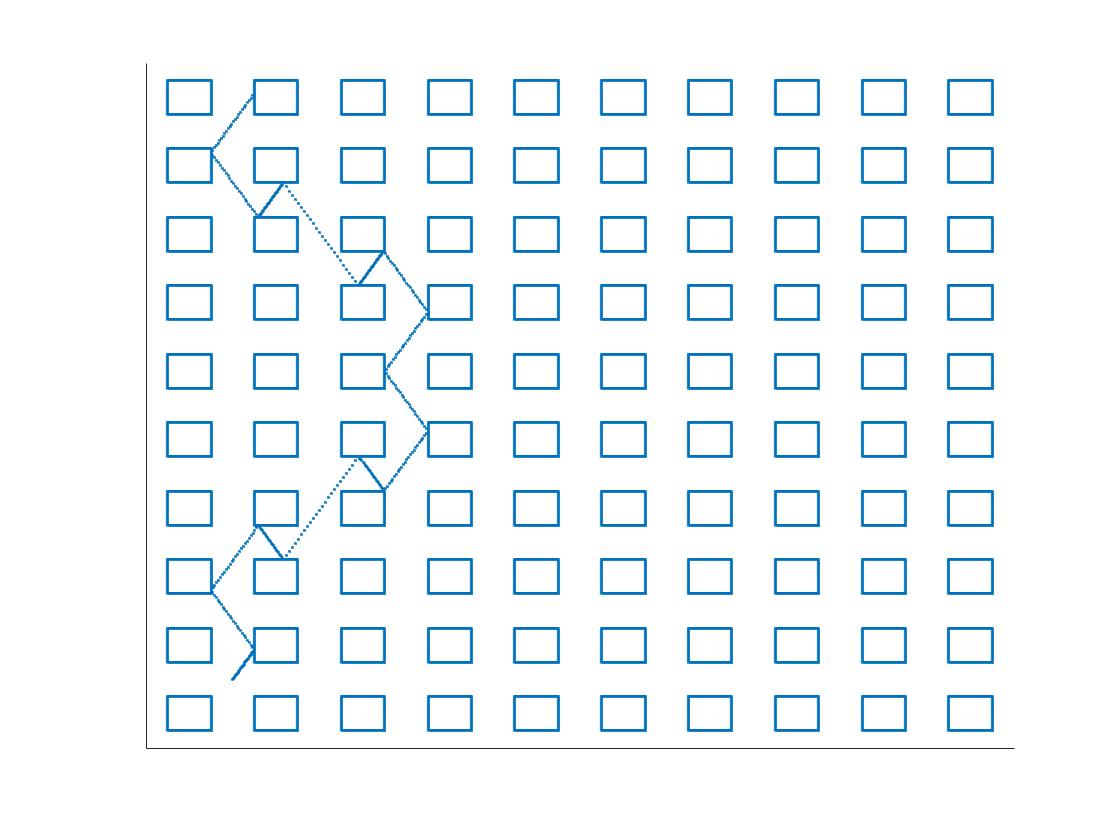}
       \caption{initial slope = 1.732}
   \end{subfigure}
   \caption{Trajectories after 15 collisions}
   \label{fig:patterns_G}
\end{figure}

 Previous authors analyze the dynamical properties of the model using techniques from algebraic topology (see \cite{Delecroix-Zorich}, \cite{Delecroix-Hubert}, \cite{Delecroix1}, \cite{Ulcigrai}).
				In this paper we take a somewhat different point of view than these works.  We assume no knowledge of algebraic topology.  Instead, we take a signal processing viewpoint and use a hidden Markov model with three hidden states to study the phenomenon we observe.

The billiard trajectory is characterized by a direction at the initial position. This direction is specified by the tangent of the angle  with a  horizontal line parallel to the plane. The slope $m$ of the initial velocity vector corresponds to the angle $\theta$, i.e.  $m = \tan \theta$. Henceforth, we can refer to either the direction $\theta$ or the initial slope $m$. \\


		\section{Background and Related Works}

			    In the periodic version of the Erhenfest
wind tree model, due to J. Hardy and J. Weber (\cite{HardyWeber}), the obstacles are identical rectangular blocks located periodically in the plane, every
obstacle centered at each point of $\mathbb{Z}^2$.  The blocks are rectangles of size $a \times b$, with $0 < a < 1, 0 < b < 1$.  We denote by $T(a,b)$ the subset of the plane obtained by removing the obstacles and name its billiard the wind tree model. 
Let $\phi_{t}^{\theta} \colon T(a, b) \rightarrow T(a, b)$ be the billiard flow: for a point $p \in T(a, b),$
the point $\phi_{t}^{\theta}(p)$ is the position of a particle after time t starting from position $p$ in direction $\theta$. 
Let $d$ be the Euclidean distance in $\mathbb{R}^2$.  The flow in direction $\theta$ is recurrent, if for almost all points $x$ in $T(a, b)$, we have $\lim \inf_{t \rightarrow \infty} d(x, \ \phi_{t}^{\theta}(x) ) = 0.$  The flow in direction $\theta$ is divergent, if for almost all points $x$ in $T(a, b)$, we have $\lim \inf_{t \rightarrow \infty} d(x, \ \phi_{t}^{\theta}(x) ) = \infty.$  Delecroix {\cite{Delecroix1}} proves the following result about the set of parameters $(a, b, \theta)$ for which the flow in $T(a,b)$ in direction $\theta$ is divergent. \\

{\bf{Theorem 1}}.    If $a$ and $b$ are either rational or quadratic of the form $1 / (1 - a) = x + y \sqrt{D}$
and $1 / (1 - b) = (1 - x) + y \sqrt{D}$
there exists a dense set $\Lambda  \in [0, 2 \pi)$ of Hausdorff dimension not smaller than $1/2$ such that for every $\theta \in \Lambda$ and every point $x$ in $T(a, b)$ with infinite forward orbit $\lim \inf_{t \rightarrow \infty} d(x, \ \phi_{t}^{\theta}(x) ) = \infty.$ In particular the flow $\phi_{t}^{\theta}$ is divergent.	\\	

 Delecroix, Hubert, and Lelièvre (\cite{Delecroix-Hubert}) determine the  rate of diffusion of the orbits which is valid for almost all direction $\theta$. \\

{\bf{Theorem 2}}. \   Then for all parameters $(a, b) \in (0,1)^2$, Lebesgue-almost all $\theta$ and every point $p$ in $T(a, b)$ (with an infinite forward orbit)
\[ \lim \sup_{T \rightarrow \infty} \frac{ \log d(p, \phi_{t}^{\theta}(p) ) }{\log T } = \frac{2}{3}. \]

 \vspace{0.15in}
	{\bf{Remark:}}  The implication is significant. 
	If we change the height and the width of the obstacle, we can get different billiard trajectories, but this does not change the diffusion rate.  To understand  the dynamical properties of the billiard motion, it is enough to take all the obstacles to be square blocks.\\

	 This theorem  is impressive because it is valid for almost all direction $\theta$ (in the sense of Lebesgue measure).  The set of rational numbers $Q$ is a set of Lebesgue measure $0$,  therefore, if $ \tan \theta  \in Q$, we cannot conclude from the  theorem whether the particle with  direction $\theta$ is recurrent or divergent. \\
	 
	   We adopt the following setting in the remainder of this manuscript:
	\begin{enumerate}
	\item All obstacles are square blocks of side length 1.
	 \item The particle starts at the origin.
         \item The centers of the blocks are the ordered pairs of all the odd integers.  
        \end{enumerate}
        This setting implies that the first obstacle to the northeast of the origin has its bottom left corner at $(0.5, 0.5)$ and its bottom right corner at $(1.5, 0.50$) with its center at $(1,1)$.  Two adjacent square blocks are at distance one apart.  There is no loss of generality in choosing this specific setting (see the {\bf{Remark}}).
       
        	\vspace{0.25in}
	
Sections \ref{demo_Part1} and \ref{demo_Part2} contain illustrations for the position of the particle after $N$ collisions with the square blocks.  Section \ref{section:Experiments} shows two experiments for  billiard trajectories after 500 collisions, from which we can identify some patterns.  These simulations suggest that there are three types of motions.
Section \ref{section:HMM} reviews some important definitions and notations in hidden Markov processes, also known as hidden Markov models.  Section \ref{section:Model_Results} describes how a 3-state hidden Markov process is adequate to model a dynamical property of the periodic wind-tree model. \\
	
	
\section{Illustrations: Part 1}\label{demo_Part1}

Initial slope = 1.414 \\
 These figures show the position of the particle after $N$ collisions with square blocks. \\

Motion of billiard  after 50 collisions, after 100, 150, and 300  collisions.

\begin{figure}[h!]
   \centering
   \begin{subfigure}[b]{0.4\linewidth}
       \includegraphics[width=\linewidth]{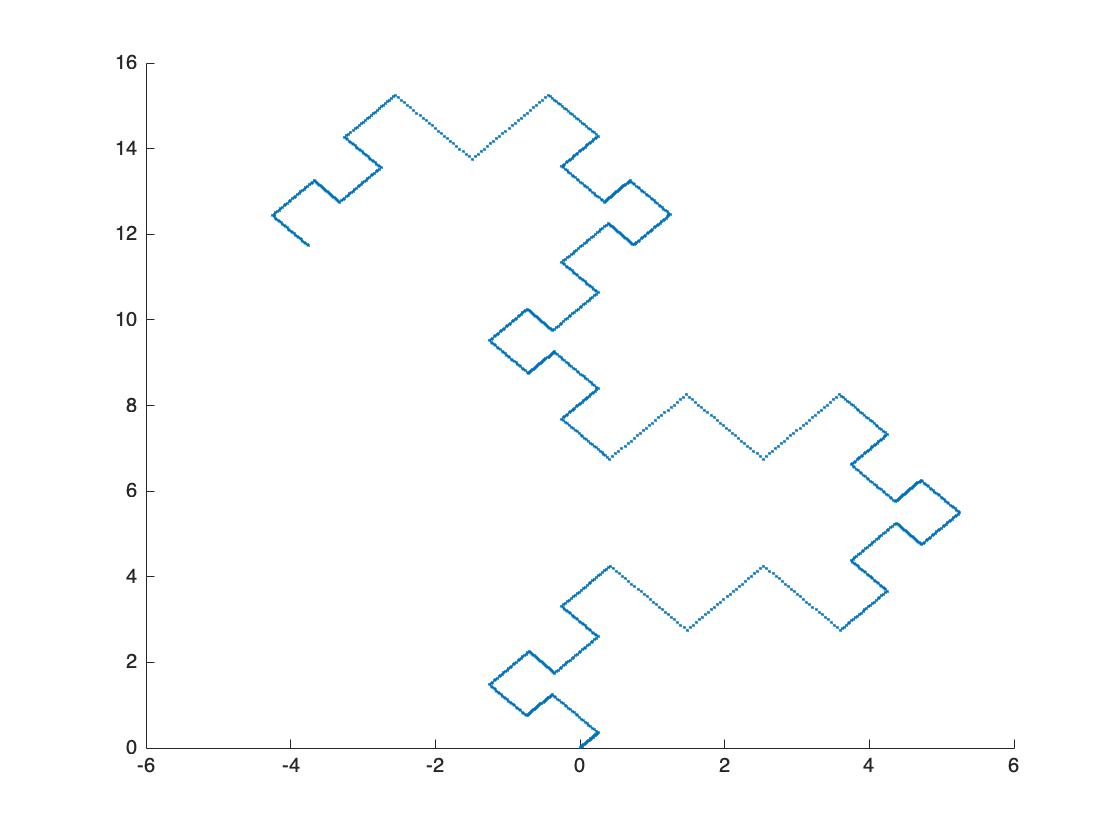}
       \caption{Pattern 1}
   \end{subfigure}
   \begin{subfigure}[b]{0.4\linewidth}
       \includegraphics[width=\linewidth]{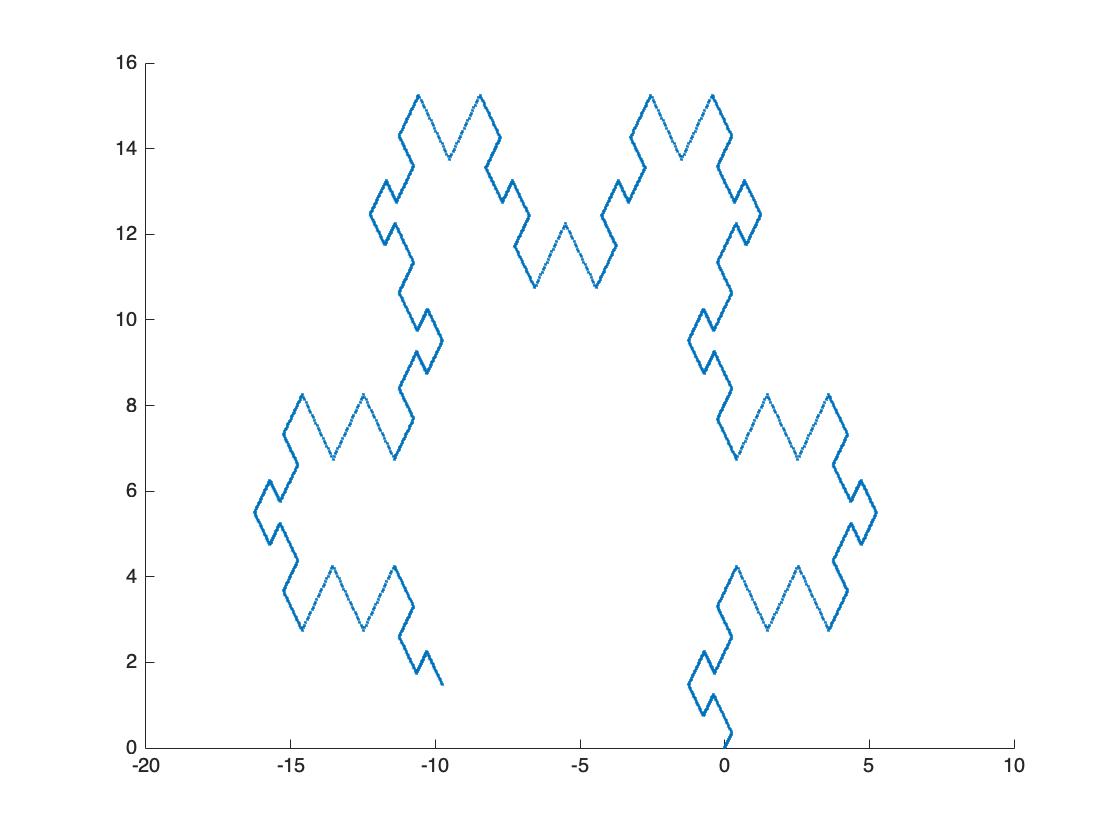}
       \caption{Pattern 2}
   \end{subfigure}
   \caption{Patterns after 50 and 100 collisions}
   \label{fig:patterns_A}
\end{figure}

\begin{figure}[h!]
   \centering
   \begin{subfigure}[b]{0.4\linewidth}
       \includegraphics[width=\linewidth]{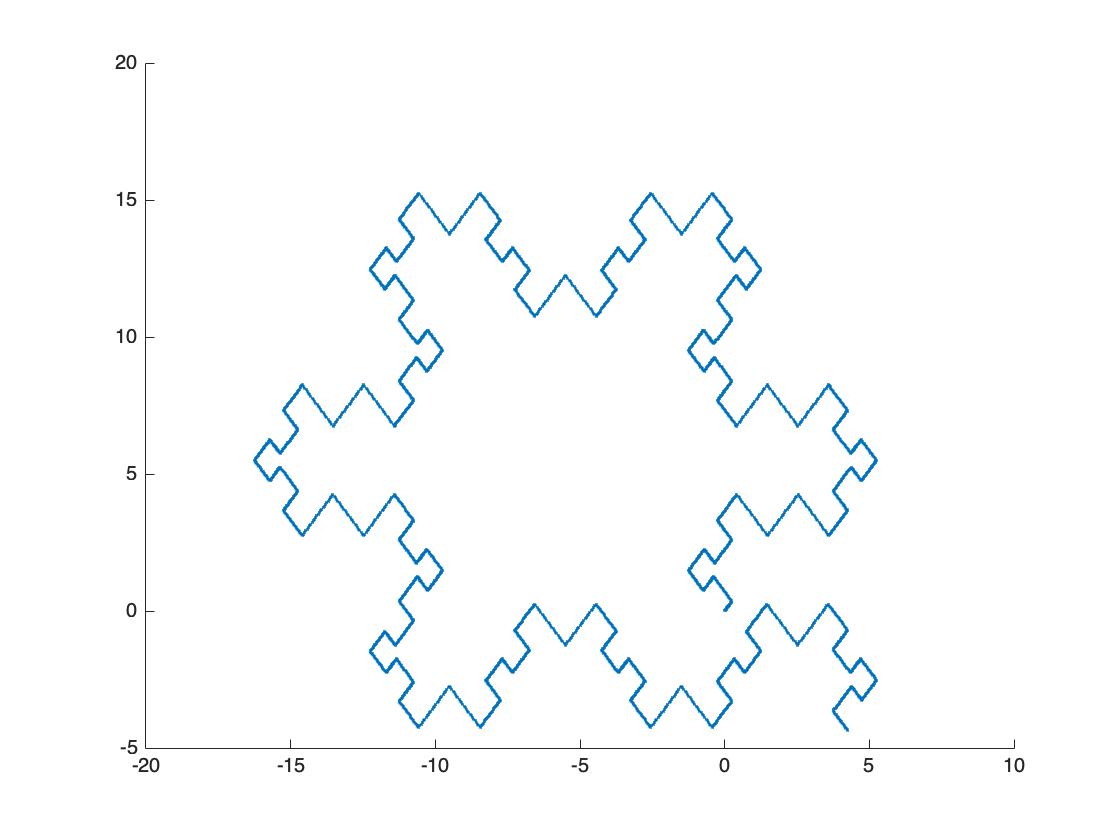}
       \caption{Pattern 3}
   \end{subfigure}
   \begin{subfigure}[b]{0.4\linewidth}
       \includegraphics[width=\linewidth]{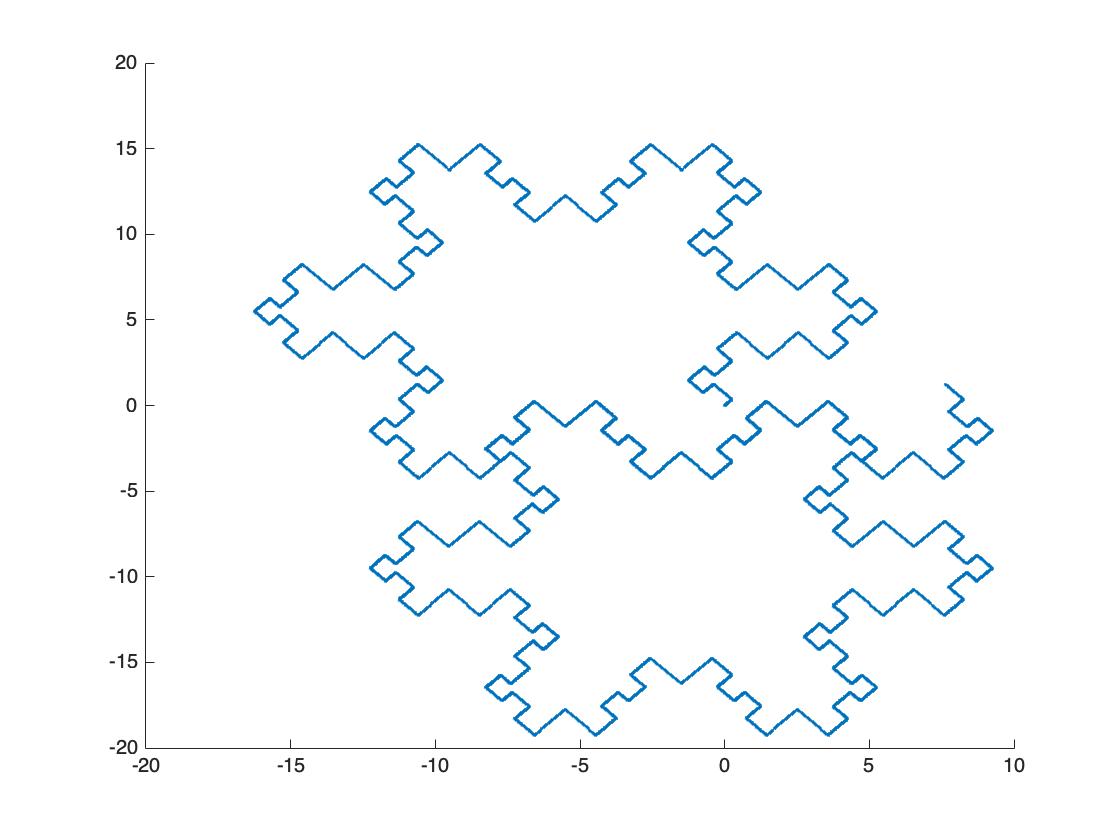}
       \caption{Pattern 4}
   \end{subfigure}
   \caption{Patterns after 150 and 300 collisions}
   \label{fig:patterns_B}
\end{figure}

\newpage

\section{Illustrations: Part 2}\label{demo_Part2}
		
Initial slope = 1.732 

\begin{figure}[h!]
   \centering
   \begin{subfigure}[b]{0.35\linewidth}
       \includegraphics[width=\linewidth]{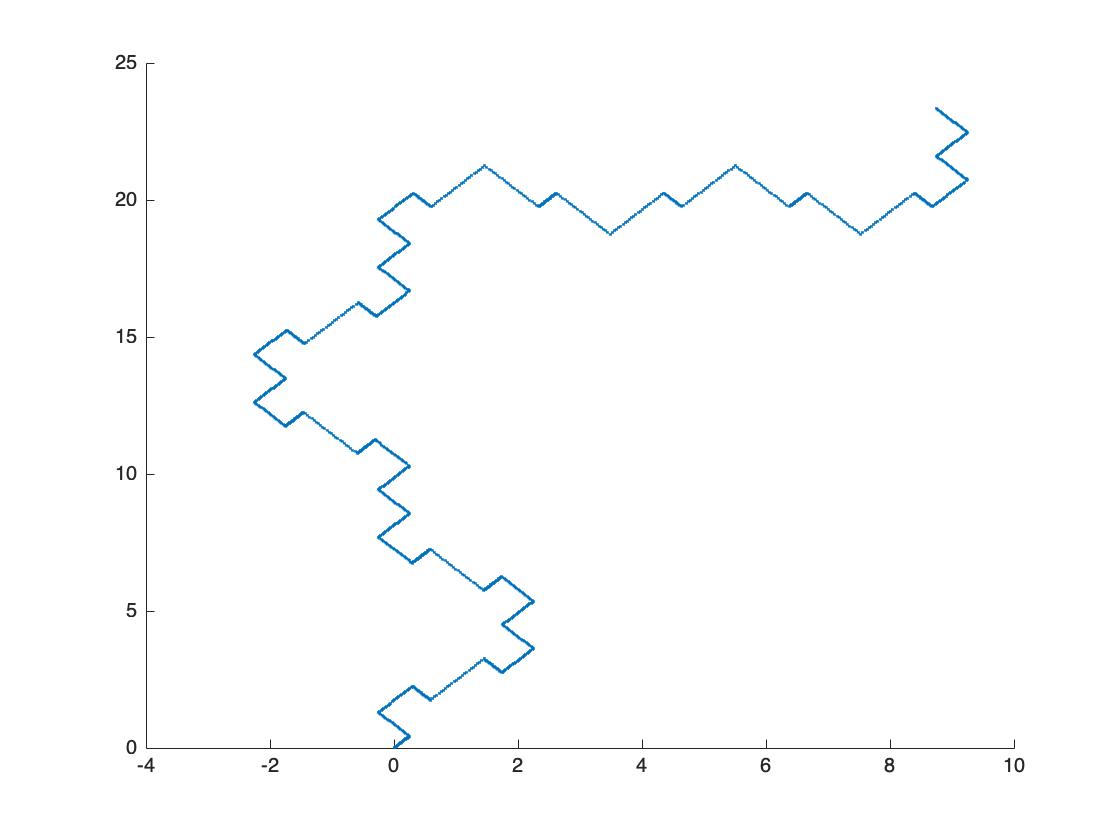}
       \caption{Pattern 1}
   \end{subfigure}
   \begin{subfigure}[b]{0.35\linewidth}
       \includegraphics[width=\linewidth]{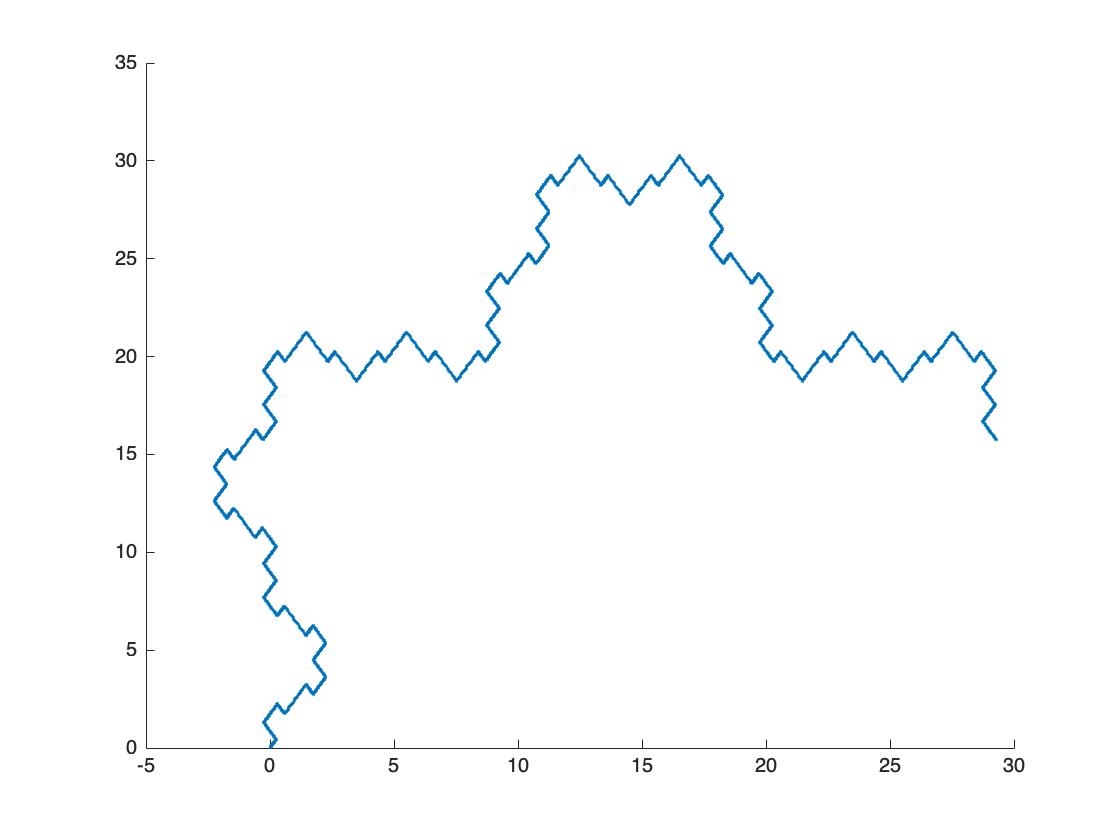}
       \caption{Pattern 2}
   \end{subfigure}
   \caption{Patterns after 50 and 100 collisions}
   \label{fig:patterns_A}
\end{figure}

\begin{figure}[h!]
   \centering
   \begin{subfigure}[b]{0.35\linewidth}
       \includegraphics[width=\linewidth]{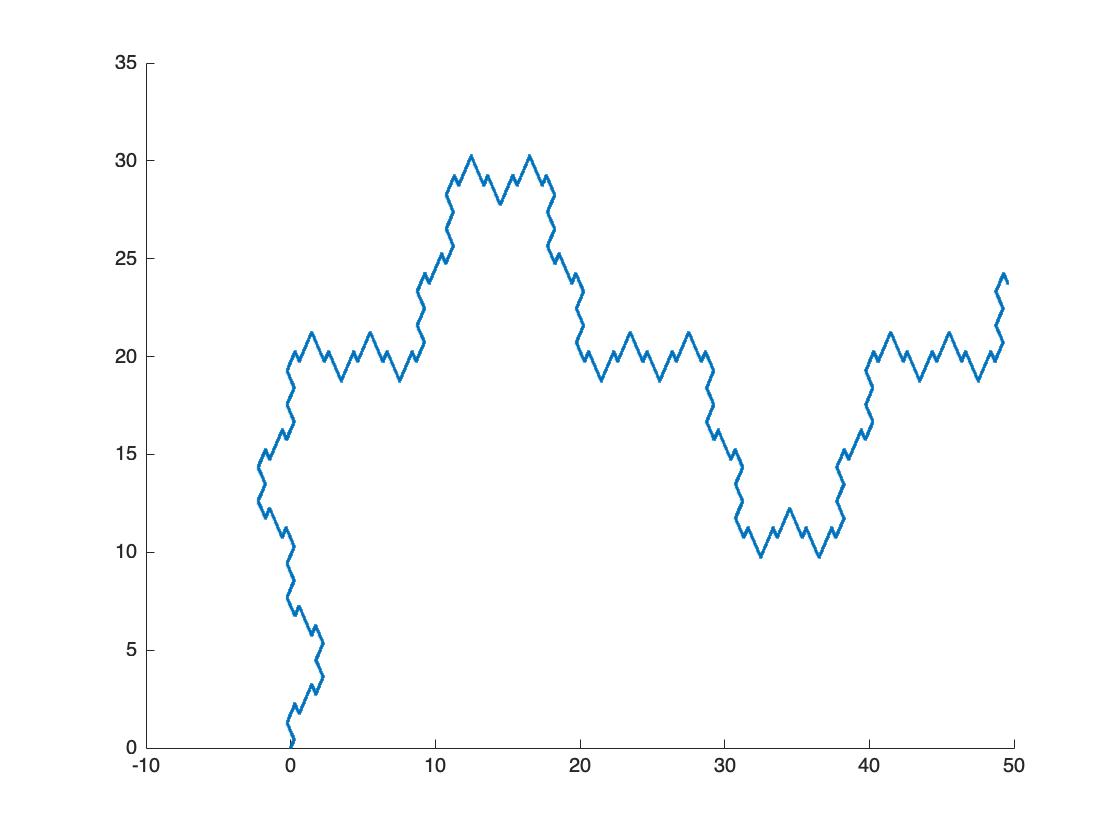}
       \caption{Pattern 3}
   \end{subfigure}
   \begin{subfigure}[b]{0.35\linewidth}
       \includegraphics[width=\linewidth]{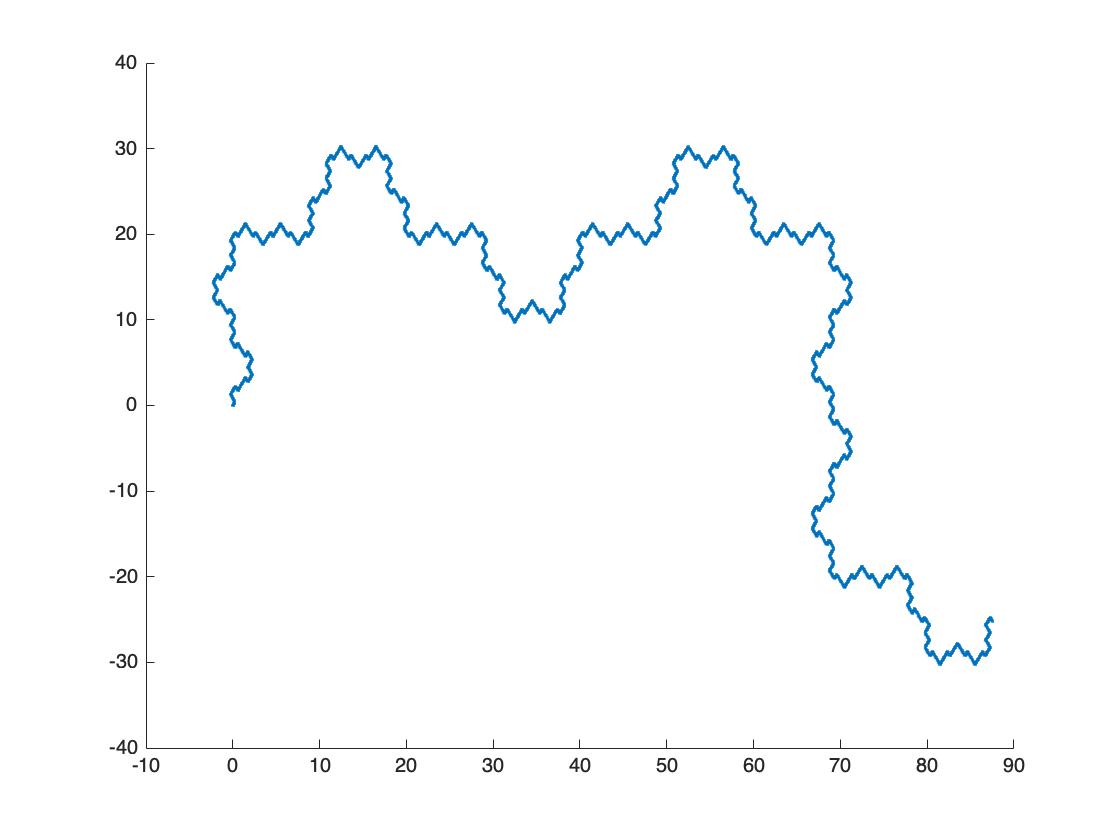}
       \caption{Pattern 4}
   \end{subfigure}
   \caption{Patterns after 150 and 300 collisions}
   \label{fig:patterns_B}
\end{figure}

 We can see that a pattern emerges after 1000 collisions.
   
 \begin{figure}[h!]
   \centering
   \begin{subfigure}[b]{0.3\linewidth}
       \includegraphics[width=\linewidth]{Pattern_4_1732.jpg}
       \caption{Pattern 4}
   \end{subfigure}
   \begin{subfigure}[b]{0.3\linewidth}
       \includegraphics[width=\linewidth]{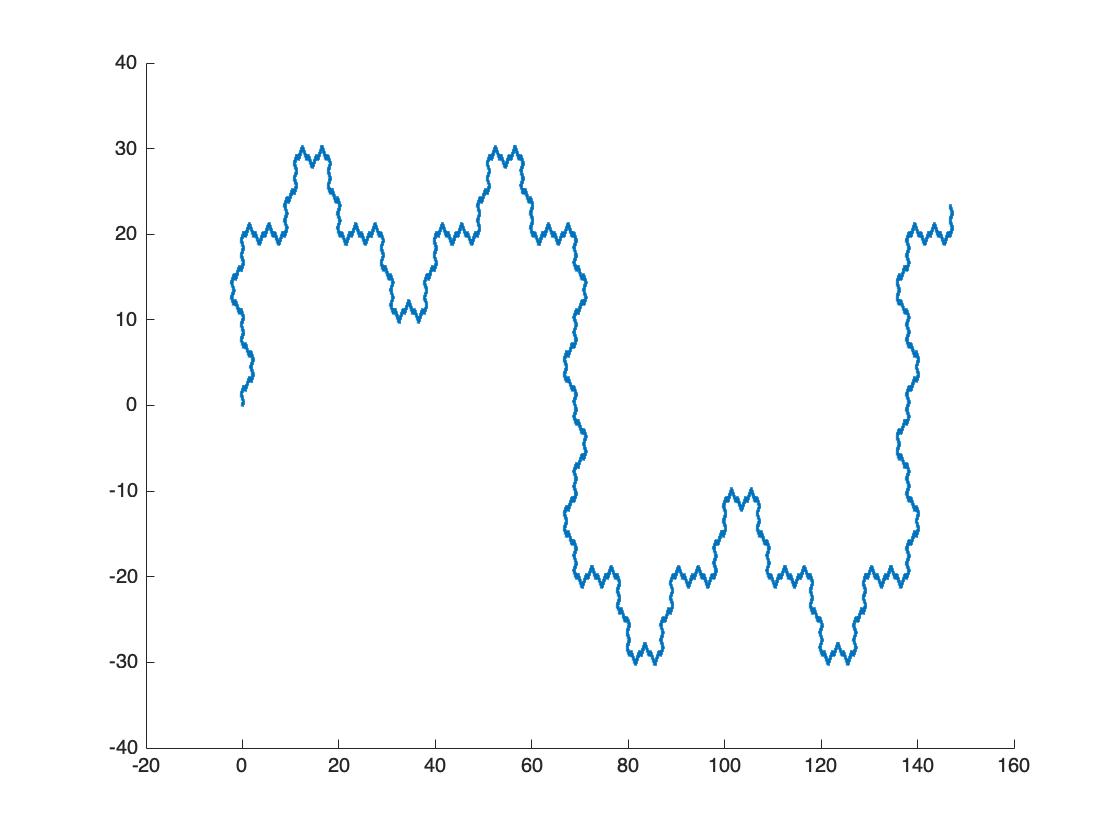}
       \caption{Pattern 5}
   \end{subfigure}
   \begin{subfigure}[b]{0.3\linewidth}
       \includegraphics[width=\linewidth]{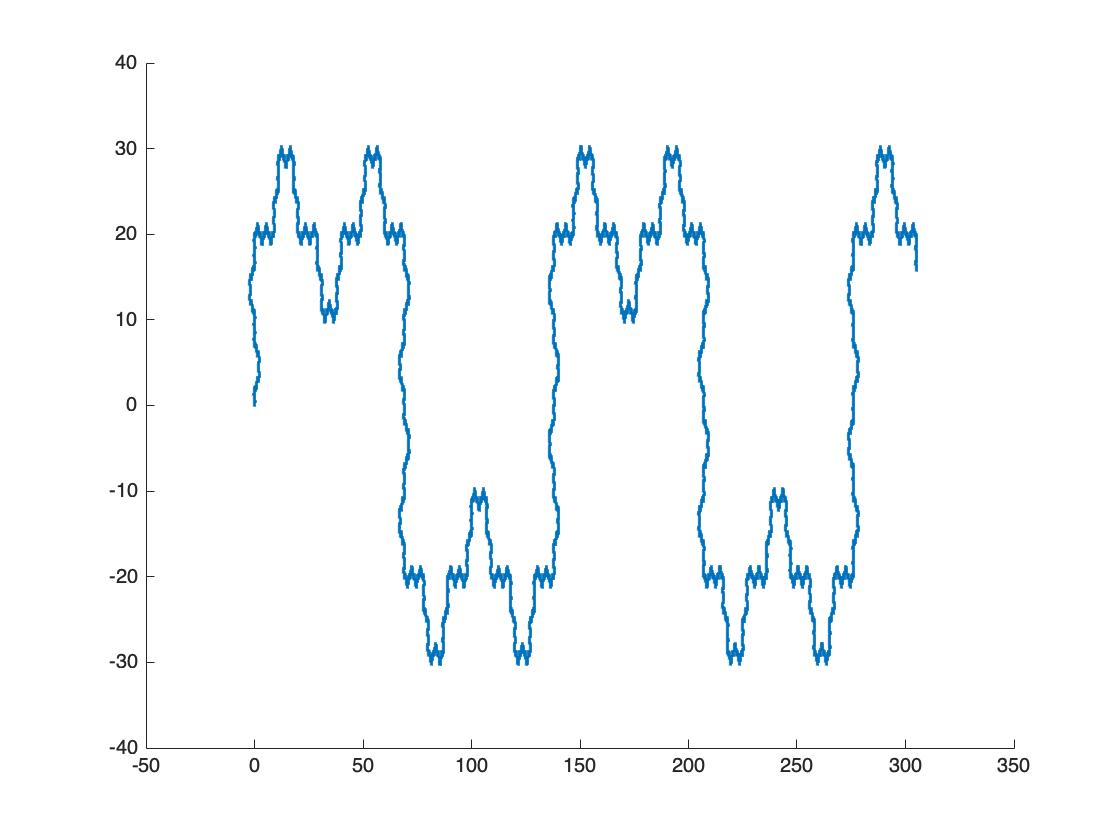}
       \caption{Pattern 6}
   \end{subfigure}
   \caption{Patterns after 300, 500, and 1000 collisions}
   \label{fig:patterns_C}
\end{figure}

\newpage

\section{Experiments: Bouncing Square Blocks}\label{section:Experiments}

Experiment 1: We observe at least 3 types of motion.\\
Billiard bouncing square blocks with given initial slope for 500 collisions. \\

(a).    Initial slope = 1.732

(b).   Initial slope = 1.618

 (c).  Initial slope = 1.414

 \begin{figure}[h!]
   \centering
   \begin{subfigure}[b]{0.3\linewidth}
       \includegraphics[width=\linewidth]{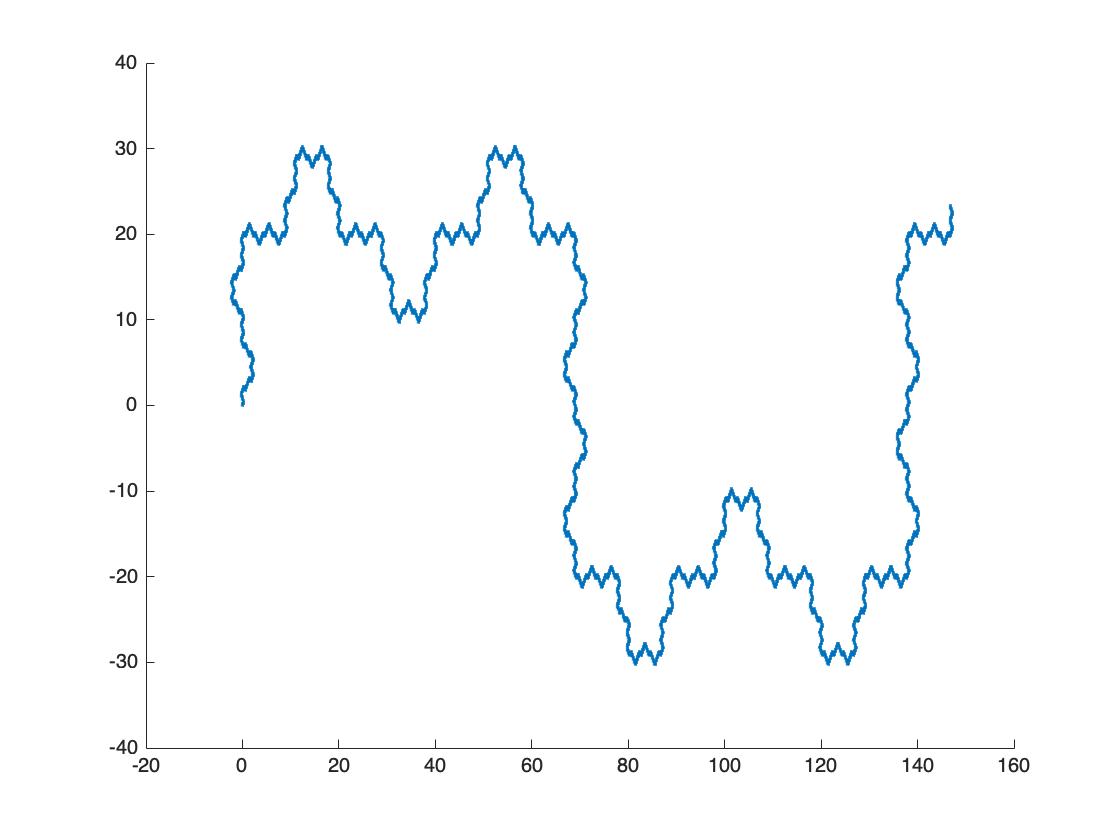}
       \caption{Initial slope = 1.732}
   \end{subfigure}
   \begin{subfigure}[b]{0.3\linewidth}
       \includegraphics[width=\linewidth]{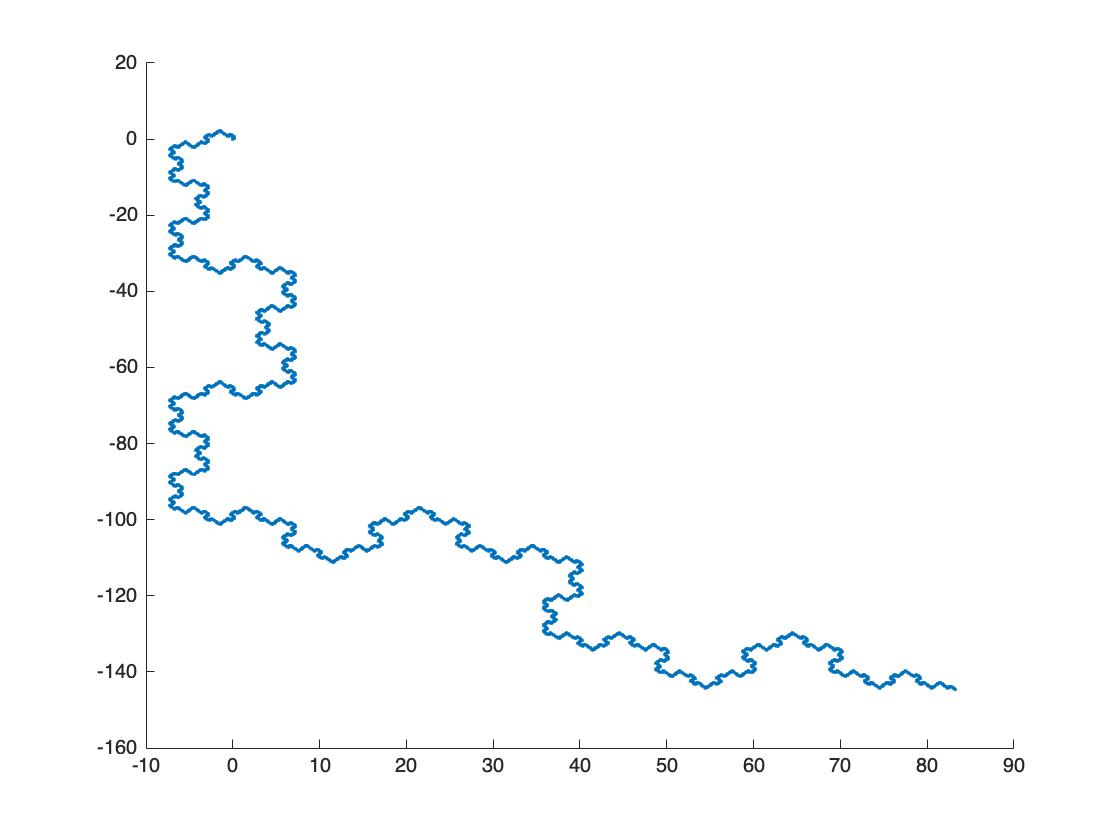}
       \caption{ Initial slope = 1.618}
   \end{subfigure}
   \begin{subfigure}[b]{0.3\linewidth}
       \includegraphics[width=\linewidth]{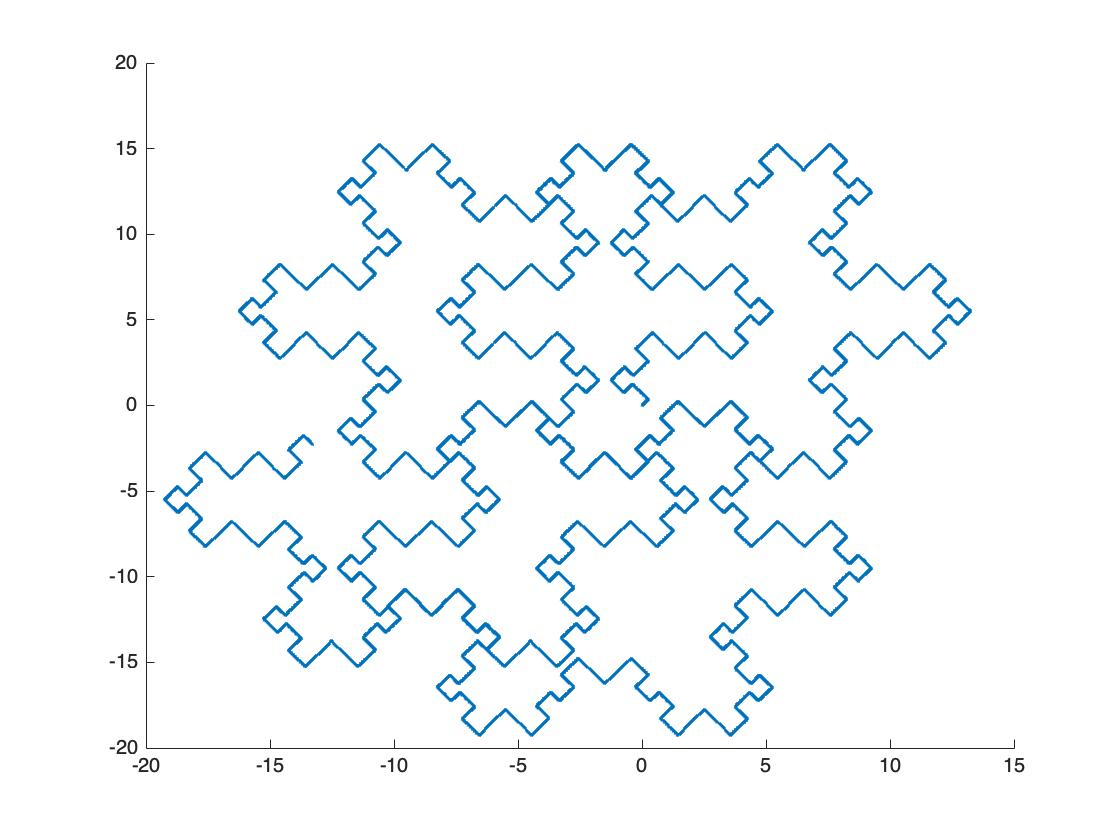}
       \caption{Initial slope = 1.414}
   \end{subfigure}
   \caption{Patterns after 500 collisions, for 3 types of motion}
   \label{fig:patterns_D}
\end{figure}
		
		We notice three types of patterns  from the first 500 collisions, as illustrated in figure 7. \\
		
 \noindent (1). Recurrent (e.g. when initial slope = 1.414) \\
		(2). Divergent in a quasi-periodic way (e.g. when initial slope = 1.732)\\
		(3). Divergent rapidly (e.g. when initial slope = 1.618)
		
		In the above classification, we use the term quasi-periodic in (2).  We say that the motion is quasi-periodic if there is a $T$, a quasi-period $\tau$, and a small number $\epsilon > 0$, such that the $y$-coordinate $y(t)$ of the particle at time $t$ satisfies
		\[ d\left(  y(t + \tau),  y(t)  \right) \leq \epsilon \]
		for all $t > T$.  The number $\tau$ can depend on the time $t$.

Experiment 2 \\
Billiard bouncing square blocks with given initial slope for 213 collisions.\\

     Initial slope = 1.718 \\
        Position at iteration 213 is   x = -0.2500, y = 0.4355 \\
        Position at first iteration is x = 0.2500, y = 0.4295 \\

\begin{figure}
			\centering
			\includegraphics[scale=0.25]{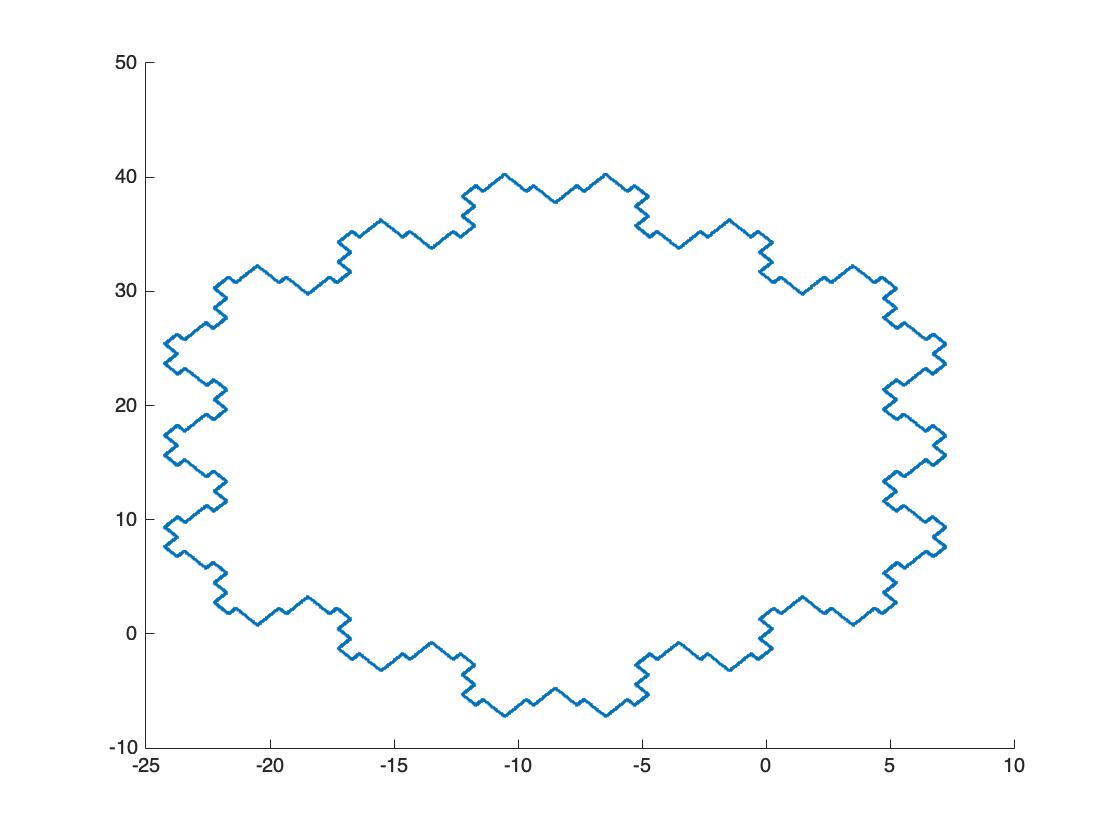}	
			\caption[ ]{Recurrent motion}
			\label{fig:bounce_recurrent}		
		\end{figure}

We observe recurrent motion.    Figure \ref{fig:bounce_recurrent} shows  
		the particle almost returns to original position.


	
	\section{Hidden Markov process}\label{section:HMM}
	Let us review some  important definitions and notations in hidden Markov processes. \\
	A sequence of discrete random variables $\{ C_t \colon t = 1, 2, 3 , \ldots \}$ is a Markov chain if  it satisfies
	the Markov property,
	\[ Pr ( \ C_{t+1} \ | C_1, C_2, C_3, \ \ldots, C_{t} ) = Pr(C_{t+1} \ | C_{t}) \]
	Consider a stochastic process $\{X(t) \colon t = 1, 2, 3, \ldots \ \}$ in which
	the probability distribution of $X(t)$ is determined by the unobserved state $C_t$ of a Markov chain.
	Let $\underline{C}^{(t)}$ and $\underline{X}^{(t)}$  be the sequence of values from time 1 to time $t$ of $C_{t}$ and $X(t)$, respectively.
	Suppose the Markov chain $\{C_t\}$ has $m$ states. Then, the stochastic process $\{X(t) \colon t = 1, 2, 3, \ldots \ \}$ is
	an $m$-state Hidden Markov model if for all $t$,
	\begin{equation}
	 Pr( X(t) \ | \underline{X}^{(t-1)}, \ \underline{C}^{(t)} ) = Pr( X(t) \ | C_t) 
	 \end{equation}
	In applications, the Marko chain $\{C_t \}$ can represent a coin process or some hidden conditions.  For example,  the daily return of stock prices $\{X(t)\}$, is influenced by  the market condition $\{C_t \}$, a Markov chain which is not directly observable.
	To illustrate the notation, we denote the history of the coin up to time $t = 4$ by $C^{(4)} = \{C_1, C_2, C_3, C_4 \}$. The main assumption for Hidden Markov model is
\[ Pr( \ X(4) = 1 \rvert \ X(1), X(2), X(3), \ C^{(4)} ) = Pr( X(4) = 1 \rvert \ C_4) \]



{\bf{Notation:}}  We denote the probability for the coin at time $t$ by
\[ \delta_i = Pr(C_t = i) \quad \mbox{ for } t = 1, 2, 3, \ldots, T.  \]
 Let $p_i$  be the probability mass function of $X(t)$ when the Markov chain is in state $i$ at time $t$.
\[ p_i(x) = Pr(X(t) = x \rvert C_t = i).  \]

As a consequence, if there are $m$ possible outcomes for the coin, then
\begin{align*} 
 Pr( \ X(t) = x \ ) & =  \sum_{i=1}^{m} Pr( X(t) = x \ \rvert \ C_t = i) Pr(C_t = i)  \\
& = \sum_{i=1}^{m} \delta_i p_i(x)
\end{align*}

Key assumption: $\{C_t\}$ is a  {\bf{homogeneous}} Markov chain. 
That means,
\[ Pr(C_{t+1} = j \rvert \ C_{t} = i) \]
does not depend on $t$.  It is a number that depends on $i$ and $j$.  We define 
\[ \Gamma(i,j) = Pr(C_{t+1} = j \rvert \ C_{t} = i). \]
$\Gamma$ is a probability transition matrix, whose entry in row $i$ and column $j$ is the probability that the coin at time $t+1$ will be in state $j$, given that the coin is in state $i$ at time $t$.

The likelihood at time $t = T$ is given by
\[ L_{T}   = Pr( \  X(1)  = x_1, X(2) = x_2, X(3) = x_3, \ \ldots \ , \  X(T) = x_{T} \ ). \]
The likelihood can be conveniently expressed in matrix notation,
\begin{equation}\label{eqn2}
L_{T} =  \vec{\delta} \ {\bf{P}}(x_1) \cdot {\bf{\Gamma}} {\bf{P}}(x_2) \cdot {\bf{\Gamma}} {\bf{P}}(x_3)   \cdot  \ \ldots \cdot {\bf{\Gamma}} {\bf{P}}(x_{T}) \ \vec{1},
\end{equation}
where ${\bf{P}}(x)$ is defined to be the diagonal matrix whose $(i, i)$ entry is $p_{i}(x)$.  Here, $\vec{1}$ denotes a column vector of ones.

To estimate the parameters of the model, we maximize the likelihood.  This can be done  by  the Baum-Welch algorithm (\cite{Baum1}, \cite{Baum2}, \cite{Rabiner}, \cite{ZucchiniMacDonald}), which is in effect one of the earliest instances of the EM algorithm (\cite{Dempster}).   In the E step of the algorithm, we infer the probabilities of the hidden states $C_{t}$ from the observations $X(t) = x_{t}$.  Then, in the M step,  we update the parameters of the model by maximizing the likelihood, conditional on the hidden states.  \\

The computation in the M step is straightforward.  The challenging part is in the E step of the EM algorithm.   We include the details in section \ref{section:hidden_states} for completeness.

\section{Model Results}\label{section:Model_Results}

In this section, we describe how a 3-state hidden Markov process is adequate to model a dynamical property of the periodic wind-tree model.

   Let $d(k)$ be the distance of the particle from the origin, after $k$ collisions with the square blocks. Let $m = \tan \theta$ be the initial slope.  Let $D(\theta)$ be the minimum distance of the particle from the origin, between the $500^{th}$ and $1000^{th}$ collision, when the initial slope is $m = \tan \theta$.     To be precise,
    \[ D(\theta) = \min\{ d(k) \colon 500 \leq k \leq 1000 \}. \]
We consider $T = 300$ values of $m$, starting from $\tan \theta_{1} = 1.4140$, goes up by increment of 0.0025, ending at $\tan \theta_{T} = 2.1615$.  Explicitly, the $300$ values are: \[1.4140, \ 1.4165, \ 1.4190, \ 1.4215, \ 1.4240, \ \ldots \ , \ 2.1590, \  2.1615. \] For each value of $m$, we compute $\log D(\theta)$.  The random process $X(t)$ is defined by \[ X(1) = \log D(\theta_{1}), X(2) = \log D(\theta_{2}), \ \ldots \ , X(T) = \log D(\theta_{T}), \]
where $\theta_1 = 1.4140, \ \theta_2 = 1.4165, \ \theta_3 = 1.4190, \ \ldots , \theta_{T} = 2.1615$.\\
To model this random process $\{X(t)\colon 1 \leq t \leq T\}$, we use a hidden Markov model with 3 states, so that conditional on state $C_{t} = j$, the random variable $X(t)$ is normally distributed with mean $\mu_j$ and variance $\sigma_{j}^2$.   The parameters of the 3-state normal-HMM are estimated from the observations $\{X(t), t = 1, 2, 3, \ldots, T\}$.

To start the EM algorithm, the probability transition matrix $\Gamma$ is initialized so that all three diagonal entries are set to $0.8$ and all the off-diagonal entries are set to $0.1$.  This initialization is arbitrary.

By applying the EM algorithm, 
 the parameters of the model after 15 iterations are: \\
 
 \begin{tabular}{|  r | l | r |} \hline 
 $\mu_1 = -0.613$ & \  $ \sigma_1 = 0.13139$   \\
$ \mu_2 = 1.9753$ & \  $ \sigma_2 = 0.10825$   \\
$ \mu_3 = 4.7825$ & \  $ \sigma_3 = 1.1217$   \\ \hline 
\end{tabular} \\

with the following $3 \times 3$ probability transition matrix $\Gamma$, \\

         0 \quad   \quad \quad 1.00     \quad    0 \\
    0.1262   \quad      0  \quad  \quad 0.8738 \\
         0  \quad \quad \quad 1.00   \quad 0 \\
         
        Notice the matrix $\Gamma$ has an appealing structure.  Only two of the nine entries are strictly between 0 and 1.

\begin{figure}[h!]
   \centering
   \begin{subfigure}[b]{0.35\linewidth}
       \includegraphics[width=\linewidth]{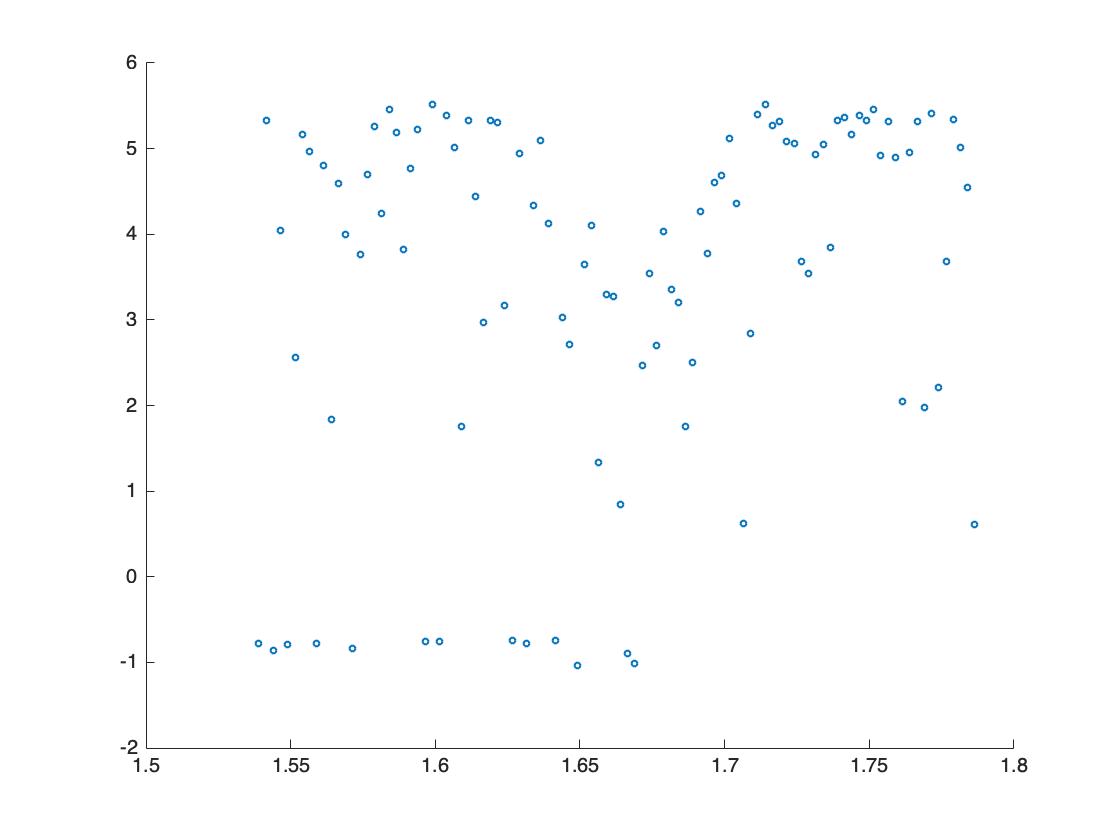}
        \caption{100 values of $X(t)$}
   \end{subfigure}
   \begin{subfigure}[b]{0.35\linewidth}
       \includegraphics[width=\linewidth]{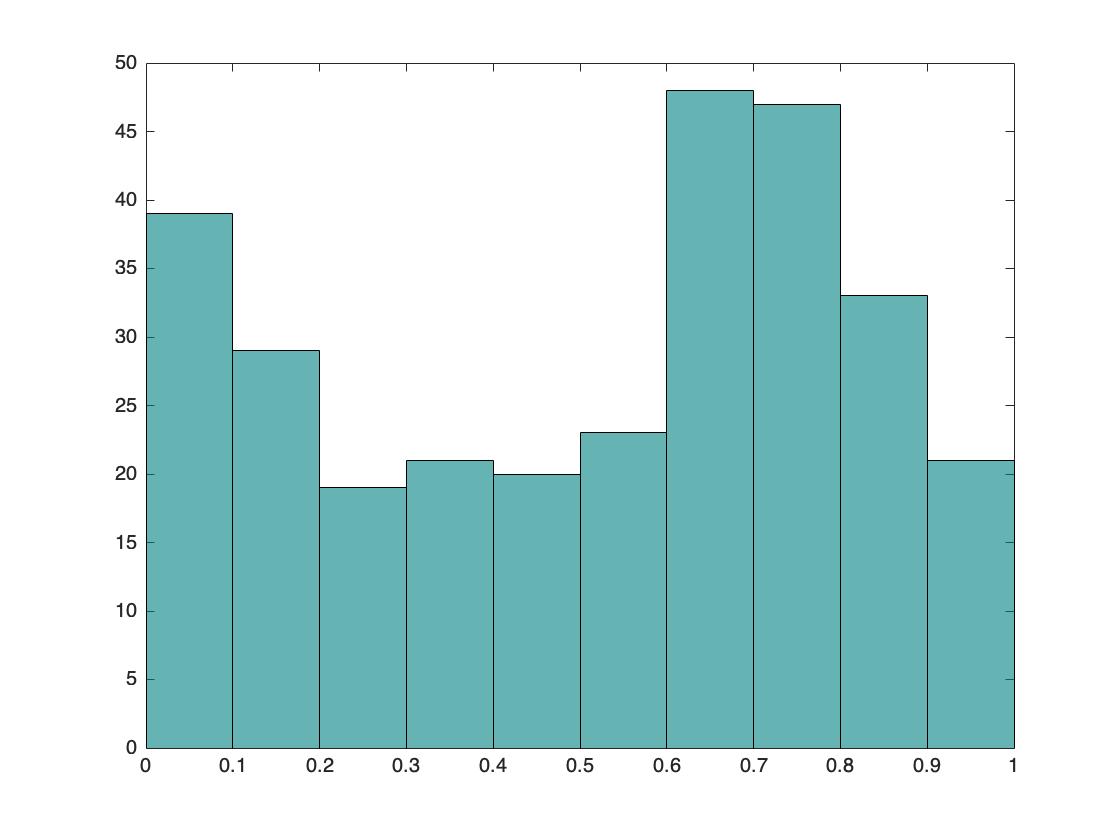}
       \caption{Plot of pseudo-residuals}
   \end{subfigure}
   \caption{The series and the pseudo-residuals}
   \label{fig:patterns_M}
\end{figure}

\vspace{0.1in}
We can interpret the states in the  3-state hidden Markov model so that the states match the 3 types of patterns we observed in Experiment 1 from section 
 \ref{section:Experiments}.
 For example, state 1 can represent recurrent motion.    We caution against the tendency to over-interpret, since we know of no convincing reason to give physical interpretation to each state.  Models need not be interpretive; they can be useful as empirical models.  We are using model as empirical model, in the sense Cox (\cite{Cox}) uses that term.\\
 
In  generalized linear models, it is common to use residuals as a tool to check the validity of the fitted model. The residuals are calculated from the expected observation and the actual observation. Pseudo-residuals are intended to fulfill the same role as the residuals, in the context of HMMs. 
The uniform pseudo-residual for the observation $x_t$ is
\[ u_t =Pr(X(t) \leq x_t)=F_{X(t)}(x_t), \]
where $X(t)$ is a continuous random variable and $F$ is its cumulative distribution function (c.d.f.). 
The diagrams of the pseudo-residuals and the series $X(t)$ are shown (for ease of viewing, the figure shows 100 values of X(t), from $t = 1.5390$ to $1.7865$).
 We check  the histogram for the pseudo-residuals, as a model diagnostic.  None of the 10 bins deviate  significantly from 30.
 The highest value among the 10 bins is 48, and the smallest  is 19.   We conclude that the model is adequate.  \\

%
%
%


\section{Probabilities of the hidden states}\label{section:hidden_states}

 In this section, we state and prove the main proposition for computing the probabilities for the hidden states $C_{t}$ in the E step of the EM algorithm.  We include the details here for completeness.
   %
   %
Before we can state the proposition, we need to describe the forward and backward probabilities.  For the sake of clear exposition, we illustrate this with 5 time periods.
    %
    %
Suppose  $T = 5$.  The likelihood at time $t = T$ is given by
  \begin{align*}
 L_{T}  & = Pr( \  X(1)  = x_1, X(2) = x_2, X(3) = x_3, X(4) = x_4, \  X(5) = x_{T} \ ) \\
           & = \vec{\delta} \ {\bf{P}}(x_1) \cdot {\bf{\Gamma}} {\bf{P}}(x_2) \cdot {\bf{\Gamma}} {\bf{P}}(x_3) \cdot {\bf{\Gamma}} {\bf{P}}(x_4) \cdot  {\bf{\Gamma}} {\bf{P}}(x_{5}) \ \vec{1},
\end{align*}
where we applied formula (\ref{eqn2}) from section \ref{section:HMM} in the last line.

{\bf{Forward probabilities}} 
 \begin{align*}
\vec{\alpha}_{1}  = \vec{\delta} {\bf{P}}(x_{1}), \quad \quad \vec{\alpha}_{2} & = \vec{\alpha}_{1} \ {\bf{\Gamma}} {\bf{P}}(x_{2}) \\
  \vec{\alpha}_{3} & = \vec{\alpha}_{2} \ {\bf{\Gamma}} {\bf{P}}(x_{3}) \\
  \vec{\alpha}_{4} & = \vec{\alpha}_{3} \ {\bf{\Gamma}} {\bf{P}}(x_{4})  
  \end{align*}
 {\bf{Backward probabilities}} 
 \begin{align*}
  \vec{\beta}_{5}  = \vec{1} , \quad \quad \vec{\beta}_{4}^{ \ T} & = \Gamma P(x_{5}) \ \vec{1}^{ \ T} \\
  \vec{\beta}_{3}^{\ T} & = \Gamma P(x_{4}) \Gamma P(x_{5}) \ \vec{1}^{ \ T} \\
  \vec{\beta}_{2}^{\ T} & = \Gamma P(x_{3}) \Gamma P(x_{4}) \Gamma P(x_{5}) \ \vec{1}^{ \ T} \\
  \vec{\beta}_{1}^{\ T} & = \Gamma P(x_{2}) \Gamma P(x_{3}) \Gamma P(x_{4}) \Gamma P(x_{5}) \ \vec{1}^{ \ T}  
\end{align*}

The forward and backward probabilities are related to the likelihood.
\begin{equation*}
 \vec{\alpha}_{t}  \vec{\beta}_{t}^{\ T} =  L_{T}   = Pr( \  X(1)  = x_1, X(2) = x_2,  \ \ldots , \ X(T) = x_{T} \ ) 
\end{equation*}

{\bf{Properties of forward and backward probabilities}} \\

The $j$-th entry of $\vec{\alpha}_{t}$ is
\[ \alpha_{t}(j) = Pr( \ X(1) = x_1, X(2) = x_2, X(3) = x_3, \ \ldots \ X(t) = x_{t} \ , \ C_{t} = j \ ) \]

The $i$-th entry of $\vec{\beta}_{t}$ is
 tells us that
 \[ \beta_{t}(i) = Pr( \ X(t+1) = x_{t+1}\ , \ \ldots \ , X(T) = x_{T}, \rvert \ C_{t} = i \ ) \]

	
    One property of the HMM is the
the conditional independence of $X(t+1)$ and $\{ X(t+2), X(t+3), \ \ldots , X(T) \}$ given $C_{t+1}$ 
\begin{align*}
 & Pr( \ X(t+1), X(t+2), X(t+3), \ \ldots , X(T) \ \rvert \ C_{t+1} \ ) \\
 = & Pr( \ X(t+1) \ \rvert \ C_{t+1} \ ) \ Pr( \ X(t+2), X(t+3), \ \ldots , X(T) \ \rvert \ C_{t+1} \ )   
 \end{align*}
 
 We will use these two properties of HMM: \\
For any integer $T \geq t + 1$, 
\begin{align*}
&  Pr( \ X(1), X(2), \ \ldots \ , X(T), \ C_{t}, \ C_{t+1} \ ) \\
= & \ Pr( \ X(1), X(2) \ \ldots \ , X(t), \ C_{t} \ ) \ Pr( \ C_{t+1} \ \rvert \ C_{t} \ ) \\
& \cdot Pr( \ X(t+1), \ \ldots \ , \ X(T) \ \rvert \ C_{t+1} \ ) \hspace{2in} (P1)
\end{align*}
and
\begin{align*}
& Pr( \ X(t+1) \ \rvert \ C_{t+1} \ ) \ Pr( \ X(t+2), \ \ldots \ , X(T) \ \rvert \ C_{t+1} \ ) \\
& = Pr( \ X(t+1), \ \ldots \ , X(T) \ \rvert \ C_{t+1} \ ) \hspace{2in} (P2)
\end{align*}

%
  {\bf{Proposition 5}} 
 \begin{align*}
  & Pr( \ C_{t-1} = j, C_{t} = k \ \rvert \ X(1), X(2), X(3), \ \ldots , X(T) \ )  \\
 = & \ \frac{ \alpha_{t-1}(j) \ \Gamma(j,k) \ p_{k}(x_{t}) \ \beta_{t}(k) }{ L_{T} } 
  \end{align*}
  
   {\underline{Proof}} of the Proposition: 
\begin{align*}
& Pr( \ C_{t-1} = j, C_{t} = k \ \rvert \ X(1), X(2), X(3), \ \ldots , X(T) \ )  \\
& = \frac{ Pr( \ X(1), \ \ldots \ , X(T), \ C_{t-1} = j, \ C_{t} = k \ ) }{ L_{T} } \\
& = Pr( \ X(1), \ \ldots \ , X(t-1), C_{t-1} = j \ ) \ Pr(C_t = k \ \rvert \ C_{t-1} = j \ ) \\
    & \cdot Pr( \ X(t), \ldots \ , X(T) \ \rvert \ C_{t} = k \ ) \cdot \frac{1}{L_{T}}, \hspace{2in} \mbox{ by } (P1) \\
 & = \alpha_{t-1}(j) \ \Gamma(j,k) \cdot Pr( \ X(t), \ \ldots , \ X(T) \ \rvert \ C_{t} = k \ ) \cdot \frac{1}{L_T} \\
 & = \alpha_{t-1}(j) \ \Gamma(j,k) \bigg(  Pr( \ X(t) \ \rvert \ C_{t} = k \ ) \ Pr( \ X(t+1), \ \ldots \ , X(T) \ \rvert \ C_{t} = k \ ) \bigg) \cdot \frac{1}{L_T} \quad \mbox{ by } (P2) \\
 & = \alpha_{t-1}(j) \ \Gamma(j,k) \ p_{k}(x_t) \ \beta_{t}(k) \cdot \frac{1}{L_T}
\end{align*}
End of Proof.

	\section{\textbf{\small Acknowledgements}} 
	
	The first author is grateful for the financial support of an Undergraduate Summer Research Program (USRP) grant from DePaul University.
	
	The second author is grateful that the Undergraduate Summer Research Program (USRP) grant from DePaul University makes this research project possible.\\
	
	
	\providecommand{\bysame}{\leavevmode\hbox to3em{\hrulefill}\thinspace}
	\providecommand{\MR}{\relax\ifhmode\unskip\space\fi MR }
	\providecommand{\MRhref}[2]{%
		\href{http://www.ams.org/mathscinet-getitem?mr=#1}{#2}
	}
	\providecommand{\href}[2]{#2}

\end{document}